\documentstyle[12pt]{article}

\let\Bbb\relax
\newfont{\Bb }{msbm10 scaled 1000}
\newfont{\Bbb}{msbm10 scaled 1200}
\newfam\eufam%
\font\euzw=eufm10 scaled 1200%
\font\euac=eufm9%
\textfont\eufam=\euzw\scriptfont\eufam=\euac%
\newcommand{\Z}{{\Bbb Z}}
\newcommand{\R}{{\Bbb R}}
\newcommand{\C}{{\Bbb C}}

\newcommand{\bR}{\mbox{{\Bbb R}}}
\newcommand{\bC}{\mbox{{\Bbb C}}}

\def\bt{\begin{thm}}
\def\et{\end{thm}}
\def\bp{\begin{prop}}
\def\ep{\end{prop}}
\def\bc{\begin{cor}}
\def\ec{\end{cor}}
\def\bl{\begin{lemma}}
\def\el{\end{lemma}}
\def\bd{\begin{dof}}
\def\ed{\end{dof}}

\def\a{\alpha}

\def\i{\iota}

\def\o{\omega}
\def\q{\theta}

\def\u{\upsilon}
\def\x{\xi}

\def\O{\Omega}

\newfont{\mcal}{eusm10 scaled \magstep1}

\def\ch{\mbox{\mcal H}}

\newfont{\goth}{eufm10 scaled \magstep1}

\def\ggl{\mbox{\goth gl}}

\chardef\tempcat=\the\catcode`\@
\catcode`\@=11
\def\cyracc{\def\u##1{\if \i##1\accent"24 i
    \else \accent"24 ##1\fi }}
\newfam\cyrfam
\font\tencyr=wncyr10 
\def\cyr{\fam\cyrfam\tencyr\cyracc}

\font\sm=cmr10
\font\sit=cmti10

\font\cmssll=cmss10 at 12 pt

\DeclareFontFamily{OT1}{msb}{}{}
\DeclareFontShape{OT1}{msb}{m}{n}
 {  <5> <6> <7> <8> <9> <10> gen * msbm
      <10.95><12><14.4><17.28><20.74><24.88>msbm10}{}
\DeclareMathAlphabet{\bubble}{OT1}{msb}{m}{n}

\def\n{\nabla}
\def\Re{\mbox {\rm Re\,}}
\def\Im{\mbox {\rm Im\,}}
\newtheorem{thm}{Theorem}
\newtheorem{prop}{Proposition}
\newtheorem{cor}{Corollary}
\newtheorem{lemma}{Lemma}
\newtheorem{dof}{Definition}
\def\scm{special complex manifold}
\def\ssm{special symplectic manifold}
\def\skm{special K\"ahler manifold}

\def\pp#1#2{{\partial #1 \over \partial #2}}

\def\wh{\widehat}

\def\ra{\rightarrow}

\def\be{\begin{equation}}
\def\ee{\end{equation}}
\def\la#1{\label{#1}} 
\def\re#1{(\ref{#1})} 
\def\arr{\begin{array}{rlll}}
\def\ea{\end{array}}
\def\bea{\begin{eqnarray}}
\def\eea{\end{eqnarray}}  

\begin{document}
\hfill J.\ Geom.\ Phys.\ (2002)

\hfill math.DG/9910091
\vskip 0.5 true cm
\begin{center}
{\LARGE Special complex manifolds}
\vskip 0.8 true cm 
{\large D.V.\ Alekseevsky\,$^1$,\  V.\ Cort\'es\,$^2$,\  
C.\ Devchand\,$^2$}
\vskip 0.5 true cm
{\small $^1$  Dept. of Mathematics, University of Hull, Cottingham Road, 
Hull  HU6 7RX, U.K.} \\
{\cyr $^{1}$  Centr ``Sofus Li'', Gen.\ Antonova 2 - 99, 117279 Moskva}\\
{\small $^2$ Mathematisches Institut der Universit\"at Bonn,
Beringstr. 1,
D-53115 Bonn\\
D.V.Alekseevsky@maths.hull.ac.uk,
vicente@math.uni-bonn.de,  
devchand@math.uni-bonn.de}
\end{center}
\vskip 0.4 true cm 
\begin{abstract} 
We introduce the notion of a special complex manifold: a complex manifold
$(M,J)$ with a flat torsionfree connection $\nabla$ such that $\nabla J$
is symmetric. A special symplectic manifold is then defined
as a special complex manifold together with a $\nabla$-parallel symplectic
form $\omega$. This generalises Freed's definition of (affine) special 
K\"ahler manifolds. We also define projective versions of all 
these geometries. 
Our main result is an extrinsic realisation 
of all  simply connected (affine or projective) 
special complex, symplectic and K\"ahler manifolds. 
We prove that the above three types of special 
geometry are completely solvable, in the sense that they are locally
defined by free holomorphic data. In fact, any special complex 
manifold is locally realised as the image of a holomorphic 1-form  
$\alpha : \bC^n \rightarrow T^*\bC^n$. Such a realisation induces a canonical
$\nabla$-parallel symplectic structure on $M$ and any special symplectic
manifold is locally obtained this way.    
Special K\"ahler manifolds are realised as complex Lagrangian submanifolds
and correspond to closed forms $\alpha$. 
Finally, we discuss the natural geometric structures 
on the cotangent bundle of a special symplectic manifold, 
which generalise the hyper-K\"ahler structure on the cotangent bundle of a
special K\"ahler manifold. 
\end{abstract} 

\medskip\noindent
{\sm 
{\sit Key words}: special geometry, special K\"ahler manifolds, 
hypercomplex manifolds, flat connections. 

\noindent   
{\sit MSC}: 53C30
\vskip .2cm  
\vfil  
\hrule width 5 truecm
\vskip .2cm

\noindent
This work was supported by the `Schwerpunktprogramm Stringtheorie' of 
the Deutsche Forschungsgemeinschaft, MPI f\"ur Mathematik (Bonn),
SFB 256 (University of Bonn), MPI  f\"ur Mathematik in den 
Naturwissenschaften (Leipzig) and 
Erwin-Schr\"odinger-Institut (Vienna).}

\newpage
\section*{Introduction}
Special K\"ahler manifolds have attracted a great deal of interest in both
string theory and differential geometry since they first arose in the
pioneering paper of de Wit and Van Proeyen \cite{dW-VP} as the allowed target 
spaces for Maxwell supermultiplets coupled to 4-dimensional $N{=}2$ 
supergravity.
These manifolds play a crucial role as admissible target spaces for scalar
and vector couplings in both {\it rigid} supersymmetric theories 
and in supergravity theories, where the supersymmetry algebra is
`locally' realised. The \skm s occurring in rigid and local supersymmetric
theories correspond respectively to the affine and projective variants
of \skm s in the mathematical literature \cite{C1,C2,F,H}. 
Special K\"ahler geometries, moreover, occur as natural geometric 
structures on certain moduli spaces.  
Projective \skm s occur, for example, as moduli spaces of Calabi--Yau
3-folds (see e.g. \cite{C1,C2}) and affine \skm s occur as moduli spaces 
of complex Lagrangian submanifolds of hyper-K\"ahler manifolds \cite{H}.    
Further, the base of any algebraic integrable system is also affine
special K\"ahler \cite{DW,F}.  

The purpose of this paper is to develop a unified perspective from 
which the various mathematical and physical approaches to special
geometry (affine versus projective, intrinsic versus extrinsic, 
 definition versus construction) can
be seen as different aspects of the same structure. 
We introduce the notion of a {\it special complex manifold} as 
a complex manifold $(M,J)$ with a flat torsionfree connection $\n$ such that 
\be 
d^{\n} J \quad =  \quad 0\, . 
\la{sc}\ee
We call it  {\it special symplectic} if, in addition, a  $\n$-parallel
symplectic form $\omega$ is specified. Further, if $\o$ is $J$-invariant, 
or equivalently, of type $(1,1)$, it is precisely a {\it special K\"ahler}
manifold in the sense of \cite{F}. More generally, we shall see 
that the Hodge components 
 $\,\omega^{11}$, $\omega^{20}$, $\omega^{02}$ of $\o$  
are closed (Proposition \ref{typeProp}). If the form $\omega^{11}$ is 
nondegenerate, it 
defines a (pseudo) K\"ahler metric $\,g=\omega^{11} \circ J\,$ on $M$
and if $\omega^{11}$ is $\nabla$-parallel (e.g., if $\omega = \omega^{11}$)  
then $(M,J,\nabla,\omega^{11})$ is a special K\"ahler manifold. 

We give an extrinsic realisation of simply 
connected special complex, symplectic and K\"ahler manifolds as immersed 
complex submanifolds of $T^* \bC^n$. 
The main property of a \scm\, used in our construction, 
is that any
affine function $f$ (i.e.\ a function satisfying $\n df =0$) can be extended
to a  holomorphic function $F$ such that $\, \Re F=f$. 
In particular, for a \ssm\ any local affine symplectic coordinate
system $(x^1,\ldots,x^n,y_1,\ldots,y_n)$ can be extended to a 
system of holomorphic functions $(z^1,\ldots,z^n,w_1,\ldots,w_n)$,
which defines a local holomorphic immersion of $M$ into $\bC^{2n}$, 
such that the special symplectic structure is induced by certain canonical 
stuctures on $\bC^{2n}$.

The fundamental example of a \scm\ $M$ is associated to a (local) holomorphic
1-form $\a=\sum F_i dz^i$ on $\bC^n$ with invertible real matrix 
$\Im \pp{F_i}{z_j}$ as follows:  
The complex manifold $M=M_\a$ is the image of the section 
$\a: \bC^n \ra T^* \bC^n = \bC^{2n}$. The flat torsionfree
connection $\n$ on $M$ is defined by the condition that the real part
$\Re F$ of any complex affine function $F$ on $\bC^{2n}$ restricts to  a
$\n$-affine function on $M$. Such a manifold $M$ carries a natural 
$\n$-parallel symplectic form $\o$ and can therefore be considered as
a \ssm\ as well. If, in addition, the 1-form $\a$ is closed (and hence locally
$\a = dF$ for a holomorphic function $F$), then $M_\a$ 
is a Lagrangian submanifold and $\o$ is of type $(1,1)$. So $M_\a$ is 
then a \skm. 
Conversely, we prove that any special complex, symplectic or K\"ahler manifold
can be locally obtained by this construction. 
More generally, we show that any totally complex holomorphic immersion $\phi$
of a complex n-manifold $M$ into $\bC^{2n}$ induces on $M$ the structure
of a \ssm. Here, we call an immersion {\it totally complex} if the intersection
$d\phi (T_pM) \cap \bR^{2n} =0$ for all $p\in M$. 
If in addition, the immersion $\phi$ is Lagrangian (i.e. $d\phi (T_pM)$ is
a Lagrangian subspace of  $T^* \bC^n $), then $M$ is a \skm.  
Our main result is that any simply connected special complex, symplectic or
K\"ahler manifold can be constructed in this fashion. In particular, 
any special K\"ahler manifold is locally defined by a holomorphic 
function $F$. This result is used  
in \cite{BC} in order to associate a parabolic affine hypersphere 
of real dimension $2n$ to any holomorphic function $F(z^1,z^2, \ldots , z^n)$
with invertible real matrix $\Im \pp{F_i}{z_j}$.  

In section 2, by including special complex manifolds   
$(M,J,\nabla)$ into a one-parameter
family   $(M,J,\nabla^\theta),\  \theta\in S^1$, we define projective
versions of special complex, symplectic and K\"ahler manifolds in terms of 
an action of $\bC^*$ on $M$ which is transitive on this family.
Our approach is based on the following observation: 
Any \scm\ $(M,J,\n )$ can be included in a one-parameter family $(M,J,\n^\q)$ 
of \scm s, with the connection $\n^\q$ defined by
\be
 \n^\q X :=   e^{\theta J}\n (e^{-\theta J}X)\, ,
\ee
where $ e^{\theta J}X = (cos \q) X + (sin \q) JX$. 
A complex manifold $(M,J)$ with a flat torsionfree connection 
$\,\n\,$ is called 
a {\it conic complex manifold~} if it admits a local holomorphic\linebreak[4]
$\mbox{\C}^*$-action  $\varphi_{\lambda}$  with differential  
$d\varphi_{\lambda}X = re^{\theta J}X$ for all $\n$-parallel
(local) vector fields $X$, where $\lambda = re^{i\theta}$.
This implies $\varphi^*_{\lambda}\n = \n^{\theta}$. 
We show that a conic complex manifold is automatically special. 

Assume that the manifold $M_\a\subset T^* \bC^n$, $\a=\sum F_i dz^i$,  
is a complex cone, i.e. it is invariant under complex scalings.
This is the case when the coefficient functions $F_i$ are homogeneous
of degree one. The induced special geometry  on $M_\a$ is then conic.  
Conversely, we prove that any conic (special) complex, symplectic
or K\"ahler manifold can be locally realised as such a cone. 
In particular, any conic \skm\  is locally described by the differential
$\a= dF $ of a holomorphic homogeneous function $F$ of degree two.  
In the simply connected case, we give a global description
of conic special manifolds in terms of holomorphic immersions.

We then define a projective special complex, symplectic
or K\"ahler manifold  as the orbit space $\overline M$ of a conic complex,
symplectic or K\"ahler manifold  $M$  by the local $\bC^*$-action, assuming
that  $\overline M$ is a (Hausdorff) manifold. 
{}From the realisation of simply connected conic manifolds as immersed
submanifolds of $T^* \bC^n$, we obtain an analogous realisation of projective
special manifolds as immersed submanifolds of complex projective space
$P(T^* \bC^n)$. From this it follows that our definition of projective
\skm s is consistent with that given by Freed \cite{F}.

Finally, we discuss the natural geometric structures 
on the cotangent bundle of a \ssm\, which are generalisations
of the known 
hyper-K\"ahler structure on the cotangent bundle of a
special K\"ahler manifold \cite{CFG,C2,F,H}. 
We prove that the cotangent bundle $N =T^*M$  of a special 
symplectic manifold  $M$ carries two  canonical complex structures : 
the standard
complex structure  $J_1$  induced by $J$  and a complex structure
$J^{\o}$, defined by $\o$ and $\nabla$.
 If the form $\o^{11}$ is 
nondegenerate, then $N =T^*M$ carries also a natural almost hyper-Hermitian
structure $(J_1, J_2 , g_N)$, i.e. a Riemannian metric $g_N$
(which is an extension of the K\"ahler metric $g = \omega^{11}
 \circ J$) and two anticommuting $g_N$-orthogonal almost complex structures 
$J_1, J_2$. This almost hyper-Hermitian structure is integrable, i.e.\
$J_1$ and $J_2$ are integrable, if and only if $\o^{11}$ is $\n$-parallel. 
In this case $(J_1, J_2 , g_N)$ is a hyper-K\"ahler structure and we
recover the known hyper-K\"ahler structure on the cotangent bundle of a
special K\"ahler manifold.
Similarly, if $\o' = \o^{20} + \o^{02}$ is nondegenerate, 
then $N =T^*M$ carries a natural almost para-hypercomplex structure
\footnote{The notion of  para-hypercomplex structure used in this paper, 
involving two commuting complex structures  and one involution  
(the product $J_1 J_2$), is a variant of the more standard notion 
consisting of two anticommuting involutions and one complex structure.}, 
that is a pair $(J_1, J_2)$ of commuting almost complex structures. 
Here $J_1$ is the standard integrable complex structure and $J_2$ is 
integrable if and only if the form $\o'$ is $\n$-parallel. 
 
We should like to thank Maximiliano Pontecorvo for
pointing out an error in our original manuscript.

\section{(Affine) special geometry} 
\subsection{Special manifolds} 
\begin{dof} A {\cmssll special complex manifold} $(M,J,\n )$ is a  
complex manifold $(M,J)$ together with a flat torsionfree connection
$\n$ (on its real tangent bundle) such that 
\[ d^{\n} J \quad =  \quad 0 \, .\]
Here the complex structure $J$ is considered as a 1-form with values in $TM$
and $d^{\n}$  denotes the covariant exterior derivative defined by $\n$.
\\
A {\cmssll special symplectic manifold} $(M,J,\n ,\omega )$ is a 
special complex manifold $(M,J,\n )$ together with a $\n$-parallel
symplectic structure $\o$. 
\\
A {\cmssll special K\"ahler manifold} is a
special symplectic manifold $(M,J,\n ,\omega )$ for which $\omega$
is $J$-invariant, i.e.\ of type $(1,1)$. The (pseudo-)K\"ahler metric
$g(\cdot , \cdot ) := \omega (J\cdot ,  \cdot)$ is called the {\cmssll 
special K\"ahler metric} of the special K\"ahler manifold 
$(M,J,\n ,\omega )$. 
\end{dof}

\noindent
{\bf Remark 1:} The evaluation of the $TM$-valued 2-form
$d^{\n} J = alt (\n J)$ on two tangent vectors $X$ and $Y$ is given by: 
\[ {d^{\n} J}(X,Y) = (\n_XJ)Y - (\n_YJ)X\, .\] 

\noindent
{\bf Remark 2:} Since we do not assume the definiteness of the metric, 
it would be more accurate to 
speak of special {\it pseudo}-K\"ahler manifolds/metrics. However, as the
signature of the metric is not relevant for our present discussion,
we shall omit the prefix {\it pseudo}. 

Given a linear connection $\n$ on a manifold $M$ and an invertible endomorphism
field $A$ on a manifold $M$, we denote by $\n^{(A)}$  the connection
defined by
\[ \n^{(A)} X = A \n (A^{-1}X)\, .\] 
Given a flat connection $\n$ on (the real tangent bundle of) a complex 
manifold $(M,J)$, we define a one-parameter family of connections
$\n^{\theta} = \n^{(e^{\theta J})}$ parametrised by the 
projective line $P^1 = \mbox{\R}/\pi \mbox{\Z}$, 
where $e^{\theta J}
= (\cos \theta ) \mbox{Id}+ (\sin \theta ) J$.  The  
connection $\n^{\theta}$ is flat, since:  
\[ \n^{\theta}X = 0 \Longleftrightarrow \n (e^{-\theta J} X) = 0 \, ,\]
where $X$ is a local vector field on $M$. 

\begin{lemma}
Let $\n$ be a flat connection with torsion $T$ on a complex manifold $(M,J)$. 
Then 
\[  \n^{\theta} = \n  + A^{\theta}\, ,\quad \mbox{where}\quad 
A^{\theta} = e^{\theta J}\n (e^{-\theta J}) = - (\sin \theta ) e^{\theta J}
\n J\, .\] 
The torsion $T^{\theta}$ of the connection $\n^{\theta}$ is given by:
\be  T^{\theta} = T + alt (A^{\theta}) 
= T -(\sin \theta )e^{\theta J}d^{\n} J\,.\la{tor} \ee
\end{lemma}

\begin{prop}\label{prop1} 
Let $\n$ be a flat torsionfree connection on a complex manifold $(M,J)$. 
Then the triple $(M,J,\n )$ defines a special complex manifold if and only
if one of the following equivalent conditions holds:
\begin{itemize}
\item[a)] $d^{\n} J = 0$. 
\item[b)] The flat connection $\n^{\theta}$ is torsionfree for some
$\theta \not\equiv 0\pmod{\pi \mbox{\Z}}$. 
\item[b')] The flat connection $\n^{\theta}$ is torsionfree for all 
$\theta$. 
\item[c)] There exists $\theta \not\equiv 0\pmod{\pi \mbox{\Z}}$ 
such that $[e^{\theta J}X, e^{\theta J}Y] = 0$ for all 
$\nabla$-parallel local vector fields $X$ and $Y$ on $M$. 
\item[c')] $[e^{\theta J}X, e^{\theta J}Y] = 0$ for all $\theta$ and 
all $\nabla$-parallel local vector fields $X$ and $Y$ on $M$. 
\item[d)] There exists $\theta \not\equiv 0\pmod{\pi \mbox{\Z}}$
such that $d (\xi \circ e^{-\theta J}) = 0$ for all 
$\nabla$-parallel local 1-forms $\xi$ on $M$. 
\item[d')] $d (\xi \circ e^{-\theta J}) = 0$ for all $\theta$ and 
all $\nabla$-parallel local 1-forms $\xi$ on $M$.  
\end{itemize}
\end{prop}

\noindent
{\bf Proof:} a) is the property defining special complex manifolds. 
Since $\n$ is torsionfree, the torsion $T^{\theta}$ of $\n^{\theta}$ is
related to $d^{\n} J$ in virtue of \re{tor} by: 
\[ T^{\theta}  = 
- (\sin \theta )e^{\theta J}d^{\n} J\, .\] 
If $\theta \not\equiv 0\pmod{\pi \mbox{\Z}}$ the endomorphism 
$(\sin\theta )e^{\theta J}$ is invertible. 
This implies the equivalence of a), b) and b'). Let $X$ and $Y$ be 
$\nabla$-parallel local vector fields. Then $e^{\theta J}X$ and  
$e^{\theta J}Y$ are $\n^{\theta}$-parallel, by the definition
of $\n^{\theta}$, and hence
\[ T^{\theta}(e^{\theta J}X,e^{\theta J}Y) 
= - [e^{\theta J}X,e^{\theta J}Y]
\, .\]  
This yields b) $\Leftrightarrow$ c) and b') $\Leftrightarrow$ c').  
For a $\nabla$-parallel local 1-form $\xi$ and $X,Y$ as above
we compute:
\begin{eqnarray*}&& d(\xi \circ e^{-\theta J})(e^{\theta J}X,e^{\theta J}Y) 
\\
 && =\   -\xi (e^{-\theta J}[e^{\theta J}X,e^{\theta J}Y]) +
e^{\theta J}X \xi (Y) - e^{\theta J}X \xi (X)\\
&& =\ 
-\xi (e^{-\theta J}[e^{\theta J}X,e^{\theta J}Y])\, 
\end{eqnarray*}
since the functions  $\xi (X)$ and $\xi (Y)$ are constant. 
This proves  the equivalences c) $\Leftrightarrow$ d) and 
c') $\Leftrightarrow$ d'), completing the proof of the proposition. $\Box$

Given a complex manifold $(M,J)$ with a flat connection $\n$, we 
say that the connection 
\[ \n^{\frac{\pi}{2}} = \n^{(e^{\frac{\pi}{2}J})} = \nabla^{(J)}  = 
\n - J\n J\] 
is its {\it conjugate connection}. 

\begin{cor}\label{equCor} 
Let $(M,J)$ be a complex manifold with a flat torsionfree connection $\n$. 
Then the following are equivalent:  
\begin{itemize}
\item[a)] $(M,J,\n)$ is a special complex manifold. 
\item[b)] The conjugate flat connection $\n^{(J)}$ is torsionfree. 
\item[c)] $[JX,JY] = 0$ for all 
$\nabla$-parallel local vector fields $X$ and $Y$ on $M$. 
\item[d)] $d (\xi \circ J) = 0$ for all 
$\nabla$-parallel local 1-forms $\xi$ on $M$. 
\end{itemize}
\end{cor}

\begin{cor} \label{scmCor} 
If $\ (M,J,\n)\ $ is a special complex manifold then $\ (M,J,\n^{\theta})\ $ 
is a special complex manifold for any $\,\theta$. If $\ (M,J,\n , \omega)\ $ 
is a special K\"ahler manifold then $\ (M,J,\n^{\theta} , \omega)\ $ is a
special  K\"ahler manifold for any $\,\theta$. 
\end{cor}

The next proposition shows that any special complex manifold also has a
canonical torsionfree {\it complex} connection, which in general is not flat. 

\begin{prop}
\label{xconnProp} Let $(M, J, \nabla )$ be a special complex manifold. Then
$D := \frac{1}{2}(\nabla + \nabla^{(J)})$ defines a torsionfree connection
such that $DJ = 0$. 
\end{prop}

\noindent
{\bf Proof:} As a convex combination of torsionfree connections, $D$ is a 
torsionfree connection. For any vector field $X$ on $M$ we compute:
\[ D_XJ = \nabla_XJ - \frac{1}{2}[J\nabla_XJ,J] = \nabla_XJ - \nabla_XJ = 0\, .
\quad\Box\] 

\begin{prop}\la{DProp}
Let $(M, J, \nabla , \omega )$ be a special K\"ahler manifold with 
special K\"ahler metric $g$ and Levi-Civita connection $\nabla^g$.
Then the following hold:
\begin{enumerate}
\item[(i)] $\nabla^g = D = \frac{1}{2}(\nabla + \nabla^{(J)})$. 
\item[(ii)] The conjugate connection $\nabla^{(J)}$ is $g$-dual to
$\nabla$, i.e.\ 
\[ Xg(Y,Z) = g(\nabla_XY,Z) + g(Y, \nabla^{(J)}_XZ)
\]  
for all vector fields $X$, $Y$ and $Z$ on $M$. 
\item[(iii)] The tensor $\nabla g$ is completely symmetric. 
\end{enumerate}
\end{prop}

\noindent
{\bf Proof:} (i) is an immediate
consequence of Proposition \ref{xconnProp} since $g = 
\omega (\cdot , J\cdot )$. (ii) follows from a direct computation,
which only uses the fact that $\omega$ is $\nabla$-parallel and $J$-invariant: 
\begin{eqnarray*} Xg(Y,Z)  & = &  X\omega (Y,JZ) 
=  \omega (\nabla_XY,JZ)  +  \omega (Y,\nabla_XJZ)\\   
& = &  g(\nabla_XY,Z) + \omega (JY,J\nabla_XJZ) 
   =  g(\nabla_XY,Z)  + g(Y,\nabla_X^{(J)}Z)\, . 
\end{eqnarray*}
Finally, to prove (iii) it is sufficient to check that $\nabla g$ is
symmetric in the first two arguments:
\begin{eqnarray*} 
& & (\nabla_Xg)(Y,Z) - (\nabla_Yg)(X,Z) \\
&& =\  Xg(Y,Z) - g(\nabla_XY,Z) - g(Y,\nabla_XZ) 
-Yg(X,Z) + g(\nabla_YX,Z)  + g(X,\nabla_YZ)\\ 
&& \stackrel{(ii)}{=}  - g(\nabla_XY,Z) + g(\nabla^{(J)}_XY,Z) +
g(\nabla_YX,Z) - g(\nabla^{(J)}_YX,Z)\\
&& =\  g(-[X,Y] +[X,Y],Z) = 0\, .\quad\Box  
\end{eqnarray*}

\begin{prop} 
\label{typeProp} 
Let $(M, J, \nabla , \omega )$ be a special symplectic manifold 
and $\o = \o^{11} + \o^{20} +\o^{02}$ the Hodge decomposition of the 
symplectic form. Then each of the components $\o^{11}$, $\o^{20}$, $\o^{02}$ 
are closed. 
\end{prop}

\noindent
{\bf Proof:} It is sufficient to check that the (1,1)-component 
$\o^{11} = \frac{1}{2}(\o + \o (J \cdot ,J\cdot ))$ is closed. 
Since $\n$ has no torsion, the exterior derivative is given by 
$d = alt \circ \n$. We compute:
\[ 2d\o^{11} = d(\o + \o (J \cdot ,J\cdot )) = 
alt \circ \n \o (J \cdot ,J\cdot ) \, .\]
Since $\n \o = 0$ for any $X_1, X_2, X_3 \in T_pM$ we obtain:
\begin{eqnarray*}  
&& 2d\o^{11}(X_1, X_2, X_3) \\
&&=\ \frac{1}{3} (\o ( (\n_{X_1}J)X_2,JX_3) +
\o ( JX_2,(\n_{X_1}J)X_3) + cycl.)\\
&&=\ \frac{1}{3} (\o ((\n_{X_1}J)X_2,JX_3) +
\o (JX_2,(\n_{X_3}J)X_1) + cycl.)\,,\
\mbox{using $d^{\n} J{=}0$\,,}\\  
&&=\  \frac{1}{3} (\o ((\n_{X_1}J)X_2,JX_3) -
\o ((\n_{X_3}J)X_1,JX_2) + cycl.) = 0\, .\qquad\Box 
\end{eqnarray*} 

\begin{prop}\label{metricProp}  
Let $(M, J, \nabla , \omega )$ be a special symplectic manifold 
and assume that $\o^{11}$ is nondegenerate. Then $(M,J,\o^{11})$ is a 
K\"ahler manifold with K\"ahler metric $g = \o^{11} (J \cdot , \cdot )$. 
$(M, J, \nabla , \o^{11})$ is special K\"ahler if and only if 
$\nabla \o^{11} = 0$. 
\end{prop}

\noindent
{\bf Proof:} It is clear that $g$ is a Hermitian metric on the complex 
manifold $(M,J)$. By the previous proposition the K\"ahler form 
$\o^{11}$ of $g$ is closed and hence $(M,J,g)$ is a K\"ahler manifold. 
The last statement is obvious. $\Box$

\subsection{Special coordinates} 
A flat torsionfree connection $\n$ on a manifold $M$ defines on it an
{\it affine structure}, i.e.\ an atlas with affine transition functions. 
A (local) function $f$ on $(M,\n )$ is called affine  
if $\n df = 0$. A local coordinate system $(x^1, \ldots , x^m)$ on $M$,
$m = \dim M$,  is 
called affine  if the $x^i$ are affine functions. 
Any affine local coordinate system $(x^1, \ldots , x^m)$ defines
a parallel local coframe $(dx^1, \ldots , dx^m)$. Conversely,
since any parallel 1-form $\alpha$ is locally the differential of an affine 
function $f$, given a parallel coframe $(\alpha^1, \ldots , \alpha^m)$ defined
on a simply connected domain $U \subset M$ there exist  
affine functions $x^i$ on $U$
such that $dx^i = \alpha^i$. The tuple  $(x^1, \ldots , x^m)$ defines
an affine local coordinate system near each point $p\in U$. 
This coordinate system is unique (as a germ, i.e.\ up to restrictions
of the coordinate domain) 
up to translations in \R$^m$. If we require in 
addition that the coordinate system is centred at $p\in U$, 
i.e.\ that $x^i(p) = 0$, then it is uniquely determined.

\begin{dof} Let $\,(M,J,\n , \omega )\,$ be a 
special symplectic manifold. A $\n$--affine local coordinate system 
$(x^1, \ldots , x^n ,y_1, \ldots ,y_n)$ on $M$ is called a 
{\cmssll\,real special coordinate system\,} if 
$\ \omega = 2\sum dx^i\wedge dy_i\,$.     
A {\cmssll\,conjugate pair of special coordinates\,} 
is a pair of holomorphic local coordinates 
$\ (z^1, \ldots , z^n)\ $ and $\ (w_1, \ldots ,w_n)\ $  such that 
$\ (x^1 = \Re z^1\,,\ldots,\, x^n = \Re z^n\,,\, y_1=\Re w_1\, ,\ldots
,\,  y_n=\Re w_n)$ is a real special coordinate system.  \end{dof}

\begin{thm} \label{scoordThm} 
\begin{enumerate}
\item[(i)] Any special symplectic manifold $(M,J,\nabla , \omega )$ 
admits a real special 
coordinate system near any point $p\in M$.  A real special 
coordinate system is unique up to an affine symplectic transformation. 
\item[(ii)] Any affine local coordinate system $(x^1, \ldots , x^n ,
y_1, \ldots ,y_n)$ on a special complex manifold 
admits a holomorphic extension, i.e.\ there exist 
holomorphic functions $z^i$ and $w_j$ with $\Re z^i = x^i$ and 
$\Re w_j = y_j$. The extension is unique up to (purely imaginary) 
translations.
\item[(iii)] Near any point of a special K\"ahler manifold there exists
a real special coordinate system which admits a holomorphic extension
to a conjugate pair of special coordinates. 
\end{enumerate} 
\end{thm}

\noindent
{\bf Proof:} The existence and uniqueness statements about
real special coordinate systems are obvious. Let 
$(x^1, \ldots , x^n ,y_1, \ldots ,y_n)$ be an affine local coordinate system
on a special complex manifold. 
Then  we 
define $\omega^i := dx^i - \sqrt{-1}J^*dx^i$. 
By Corollary \ref{equCor} $J^*dx^i = dx^i \circ J$ is closed. This implies
that $\omega^i$ are closed 1-forms of type $(1,0)$ and are hence closed
holomorphic 1-forms. 
So there exist local holomorphic functions $z^i$ 
such that $\omega^i = dz^i$.  By adding real constants we can arrange that
$\Re z^i = x^i$. Similarly, there exist local holomorphic functions 
$w_j$ such that $\Re w_j = y_j$.  The uniqueness  
statement concerning this holomorphic extension is obvious. 
We claim that in the case of special K\"ahler manifolds,
real special coordinates can be chosen such 
that the $dz^i$, as well as the $w_j$, are linearly  
independent (over $\mbox{\C}$). To see this, let us 
first observe that the $dx^i$ and $dy_j$ define a Lagrangian 
splitting of $T_p^*M$ with respect to $\omega^{-1}$ for any point $p$ 
in the coordinate domain: $T_p^*M = L_x \oplus L_y$, 
where $L_x = \mbox{span} \left\{ dx^1, \ldots , dx^n\right\}$ and 
$L_y =  \mbox{span} \left\{ dy_1, \ldots , dy_n\right\}$. 
The $J$-invariance of the symplectic (K\"ahler) form $\omega$ 
implies  the existence of a Lagrangian splitting of the form  
$T^*_pM = L \oplus J^*L$. Since any two Lagrangian splittings of a symplectic
vector space are related by a linear symplectic transformation, this 
shows that the real special coordinates $x^1,\ldots y_n$ near $p$
can be chosen such that the corresponding Lagrangian subspaces
$L_x$,  $L_y$ satisfy $L_x \cap J^*L_x = L_y \cap J^*L_y = 0$ 
at the point $p$ and hence on a coordinate domain containing $p$. 
The equation $L_x \cap J^*L_x = 0$ forces the  
$dz^i = dx^i -\sqrt{-1}J^*dx^i$ 
to be linearly independent. So the $z^i$ define local holomorphic 
coordinates on the special K\"ahler manifold.  Similarly, 
as a consequence of the equation $L_y\cap J^*L_y = 0$, the   
$w_j$ are local holomorphic coordinates.  
$\Box$

\subsection{The extrinsic construction of special manifolds} 
\label{extrsmSec} 
As in \cite{C2}, we consider the following fundamental algebraic data: 
the complex vector space $V = T^* \mbox{\C}^n = \mbox{\C}^{2n}$ with 
canonical coordinates 
$(z^1,\ldots ,z^n,w_1,\ldots , w_n)$ and standard complex symplectic form 
$\Omega =  \sum_{i=1}^n dz^i \wedge dw_i$, the standard real structure 
$\tau: V \rightarrow V$ with fixed point set $V^{\tau}= T^* \mbox{\R}^n$. Then 
$\gamma:= { \sqrt{-1}} \; \Omega(\cdot ,\tau \cdot)$ defines a Hermitian
form of (complex) signature $(n,n)$. 

Let $M$ be a connected complex $n$-fold. A holomorphic immersion 
$\phi: M \rightarrow V$  is 
called {\it nondegenerate} (respectively, {\it Lagrangian}) if 
$\phi^* \gamma$ is nondegenerate (respectively, if $\phi^* \Omega = 0$).  
If $\phi$ is nondegenerate, then $\phi^*\gamma$  defines a, 
possibly indefinite,  K\" ahler 
metric $g$ ($=\Re \phi^*\gamma$) on the complex manifold $M$. 
The correponding K\"ahler form 
$g(\cdot ,J\cdot )$ is a $J$-invariant symplectic form on $M$, 
where $J$ denotes the complex structure of $M$.  $\phi$  
is called {\it totally complex} if $V^{\tau} 
\cap d\phi T_pM = 0$ for all $p\in M$. 
\begin{lemma}\label{immersionLemmma}
A holomorphic immersion 
$\phi: M \rightarrow V$  is totally complex if and only if its real part
$\Re \phi: M \rightarrow V^{\tau}$ is an immersion. 
\end{lemma}

\noindent
{\bf Proof:} Let $\phi: M \rightarrow V$ be a totally complex holomorphic
immersion. Restricting, i.e. pulling back via $\phi$,  the functions 
$x^i := \Re z^i$ and $y_j := \Re w_j$ to $M$ 
we obtain $2n$ functions on $M$ 
with everywhere linearly independent differentials. In fact, let 
$\alpha = \sum a_idx^i + \sum b^jy_j$ be a real linear combination which 
vanishes on the complex $n$-dimensional linear subspace 
$d\phi T_pM \subset V$. Then, since
$\alpha$ is real, it must also vanish on $\tau d\phi T_pM$. Now we can
conclude that $\alpha = 0$ since, by our assumption on $\phi$, 
$d\phi T_pM \cap \tau d\phi T_pM = V^{\tau}\cap d\phi T_pM \oplus 
iV^{\tau}\cap d\phi T_pM = 0$ and, therefore, $V = d\phi T_pM 
\oplus \tau d\phi T_pM$.  This shows that the functions $x^i$ and $y_j$
restrict to local coordinates on $M$ and, hence, that $\Re \phi$ is an
immersion. Conversely, let $\phi: M \rightarrow V$ be a holomorphic
immersion such that $\Re \phi: M \rightarrow V^{\tau}$ is an immersion. 
We have to show that $V^{\tau} 
\cap d\phi T_pM = 0$ for all $p\in M$. Suppose, that $X \in T_pM$
and $d\phi X \in V^{\tau}$. Then we have that $0 = \Im d\phi X =
- \Re \sqrt{-1}d\phi X = - \Re d\phi JX$. This implies that
$JX = 0$, because $d\Re \phi = \Re d\phi$ is injective. This shows that
$X = 0$ proving  $V^{\tau} \cap d\phi T_pM = 0$. $\Box$

A holomorphic totally complex immersion $\phi$  induces a flat 
torsionfree connection on the real tangent bundle of $M$ as follows. 
Since $\Re \phi$ is an immersion, by Lemma \ref{immersionLemmma}, restricting 
the functions $x^i = \Re z^i$ and $y_j = \Re w_j$ to $M$ we obtain 
local coordinates, which induce a flat torsionfree connection
$\nabla$ on $M$. Moreover, $2\sum dx^i\wedge dy_i$ restricts to a 
$\nabla$-parallel symplectic form $\omega$ on $M$.  We call
$\nabla$ and $\omega$ the induced connection and the induced symplectic form 
respectively. Now we can easily prove: 

\begin{thm}\label{ssconstrThm}
Let $\phi$ be a totally complex holomorphic immersion of a complex
manifold $(M,J)$ into $V = T^* \mbox{\C}^n$, $n = \dim_{\mbox{\C}} M$,  
$\nabla$ the induced connection and $\omega = 2\phi^*(\sum dx^i\wedge dy_i)$ 
the induced symplectic form. Then the following hold:
\begin{enumerate}
\item[(i)] $(M,J,\nabla , \omega )$ is a special 
symplectic manifold. 
\item[(ii)] The pull back via $\phi$ of the functions 
$(x^1 = \Re z^1, \ldots , x^n = \Re z^n, y_1 = \Re w_1, \ldots ,$\\ 
$y_n = \Re w_n)$ of $V$ defines a real special coordinate system 
around each point of $M$.  
\end{enumerate}
\end{thm}

\noindent
{\bf Proof:} We have to prove that $d^{\n}J = 0$. By Corollary \ref{equCor}, 
it is sufficient to check that the 1-forms $dx^i \circ J$ and $dy_j \circ J$
are closed. This follows immediately from the fact that the 1-forms 
$dz^i = dx^i - \sqrt{-1}dx^i\circ J$ and $dw_j = dy_j - \sqrt{-1}dy_j\circ J$ 
are closed. $\Box$ 

\noindent
The next proposition clarifies the relation between the three notions defined
above. 

\begin{prop}\label{2->3Prop} 
Let $\phi$ be a holomorphic immersion of a complex $n$-fold 
$M$ into $V = T^* \mbox{\C}^n$. The following conditions are
equivalent:  
\begin{enumerate}
\item[(i)] $\phi$ is Lagrangian and nondegenerate. 
\item[(ii)] $\phi$ is Lagrangian and totally complex. 
\end{enumerate}
\end{prop}

\begin{thm}\label{sKconstrThm}
Let $\phi$ be a  holomorphic nondegenerate Lagrangian immersion of 
a complex manifold $(M,J)$ into $V$ inducing the K\"ahler metric $g$ on $M$.  
The immersion $\phi$ is totally complex and hence induces 
also the data $(\nabla , \omega )$ on $M$. Moreover, the following hold:
\begin{enumerate}
\item[(i)] $(M,J,\nabla , \omega )$ is a special K\"ahler manifold.
\item[(ii)] $\omega$ coincides with the K\"ahler form of $g$, i.e.\ 
$\omega = g (\cdot , J\cdot )$.
\item[(iii)] The pull back via $\phi$ of the canonical  
coordinates  $(z^1, \ldots , z^n, w_1, \ldots , w_n)$ of $V$ defines 
a conjugate pair of special coordinates around each point of $M$. 
\end{enumerate}
\end{thm}

\noindent
{\bf Proof:} Thanks to Proposition \ref{2->3Prop} and Theorem \ref{ssconstrThm}
it is sufficient to prove that $g(\cdot , J\cdot ) = \omega = 2\phi^*
(\sum dx^i\wedge dy_i)$. A straightforward computation, which only uses
the definition of $g$, shows that 
\begin{equation} 
\label{gEqu} 2g(\cdot , J\cdot ) = \omega + J^*\omega\, .\end{equation} 
On the other hand, since $\phi$ is Lagrangian, we know also that
\[ 0 = 2 \Re \phi^*\Omega = \omega - J^*\omega\, .\] 
This  implies that $g(\cdot , J\cdot ) = \omega$. $\Box$  

Now we will show that any simply connected special  (complex, symplectic or 
K\"ahler) manifold
arises by the construction of Theorem \ref{ssconstrThm} or 
Theorem \ref{sKconstrThm}. 
\begin{thm}\label{ssfundThm}
\begin{enumerate}
\item[(i)] Let  $(M,J,\nabla)$ be a simply connected special 
complex manifold of complex dimension $n$. Then there exists a 
holomorphic totally complex  immersion $\phi : M \rightarrow 
V = T^* \mbox{\C}^n$ inducing the connection $\n$ on $M$. 
Moreover, $\phi$ is unique up to
an affine transformation of $V$ preserving the real structure  $\tau$. 
Here the real structure is considered as a (constant) field of antilinear 
involutions on the tangent spaces of $V$. Finally, $\omega = 2\phi^*
(\sum dx^i\wedge dy_i)$ is a $\n$-parallel symplectic structure 
defining on $(M,J,\n )$ the structure of special symplectic manifold.  
\item[(ii)] Let  $(M,J,\nabla , \omega )$ be a simply connected special 
symplectic manifold of complex dimension $n$. Then there exists a 
holomorphic totally complex  immersion $\phi : M \rightarrow 
V = T^* \mbox{\C}^n$ inducing the connection $\n$ and the
symplectic form $\omega$ on $M$. Moreover, $\phi$ is unique up to
an affine transformation of $V$ preserving the complex symplectic 
form $\Omega$ and the real structure  $\tau$.  
\item[(iii)]
Let  $(M,J,\nabla , \omega )$ be a simply connected special 
K\"ahler manifold of complex dimension $n$ then there exists a 
holomorphic nondegenerate Lagrangian (and hence totally complex) 
immersion $\phi : M \rightarrow 
V = T^* \mbox{\C}^n$ inducing the K\"ahler metric $g$, the 
connection $\n$ and the symplectic form $\omega = 
2\phi^*(\sum dx^i\wedge dy_i) = g(\cdot , J\cdot )$ on $M$. Moreover 
$\phi$ is unique up to 
an affine transformation of $V$ preserving the complex symplectic 
form $\Omega$ and the real structure  $\tau$. Here the real structure
is considered as a field of antilinear involutions on the tangent spaces
of $V$.   
\end{enumerate}
\end{thm}

\noindent
{\bf Proof:} We prove (ii) and (iii). The proof of (i) is similar.  
By Theorem \ref{scoordThm} there exist real special coordinates
near each point of $M$. Since $M$ is simply connected, we can choose
these local coordinates in a compatible way obtaining globally
defined functions $x^i$ and $y_j$ on $M$ such that
$(x^1,\ldots , x^n, y_1, \ldots ,y_n)$ is a real special coordinate
system near each point of $M$. Then again by Theorem \ref{scoordThm}
and the simple connectedness of $M$ we can holomorphically 
extend these functions, i.e.\ there exist globally defined holomorphic 
functions $z^i$ and $w_j$ such that $\Re z^i = x^i$ and $\Re w_j = y_j$.
Moreover, if $(M,J,\nabla , \omega )$ is special K\"ahler we can
assume that $(z^1,\ldots ,z^n, w_1, \ldots ,w_n)$ form a conjugate pair
of special coordinates. 
We define the holomorphic map  
\[ \phi := (z^1,\ldots ,z^n, w_1, \ldots ,w_n) :M \rightarrow  
\mbox{\C}^{2n} = V\, .\]
The fact that $\phi$ is a totally complex immersion follows from the 
linear independence of $(dx^1,\ldots , dx^n, dy_1,\ldots ,dy_n)$.
This proves the existence statement in (ii).  
To prove (iii) we need to check that $\phi$ is Lagrangian, i.e. that
the holomorphic 2-form $\O := \sum dz^i \wedge dw_i = 0$.
This follows from the $J$-invariance of $\o=2\sum dx^i\wedge dy_i$, since
$2\,\Re \O = \o - J^*\o\,$ and 
$\,2\Im \O = J\cdot\o = 2\sum Jdx^i \wedge dy_i + 2\sum Jdy_i \wedge dx^i$. Here
the $\cdot$ stands for the natural action of $\ggl(E)$ on $\wedge^2 E^*$,
where $E= T_p M,\, p\in M$.   
The uniqueness statement is a consequence of the uniqueness statement
in Theorem \ref{scoordThm}. $\Box$

We will call a holomorphic $1$-form $\sum F_idz^i$ on an 
open subset  $U \subset \mbox{\C}^n$ {\it regular} if the real matrix  
$\Im \partial F_i/\partial z^j$ is invertible. A holomorphic 
function $F$ on $U$ is called {\it nondegenerate} if its differential 
$dF$ is a regular holomorphic $1$-form. 
Any holomorphic $1$-form $\phi$ on a domain $U \subset \mbox{\C}^n$
can be considered as a holomorphic immersion 
\[ \phi : U \rightarrow V = T^*\mbox{\C}^n\, .\]
So it makes sense to speak of totally complex or Lagrangian
holomorphic $1$-forms. 
\begin{lemma}
\label{reg->Lemma} 
Let $\phi$ be a holomorphic $1$-form. Then the following hold: 
\begin{enumerate}
\item[(i)] $\phi$ is totally complex if and only if it is regular. 
\item[(ii)] $\phi$ is Lagrangian if and only if it is closed. 
\end{enumerate}
\end{lemma}

\noindent
{\bf Proof:} (ii) is a well known fact from classical mechanics. 
To see (i) let $\phi = \sum F_idz^i$ be a holomorphic $1$-form
on a domain $U \subset \mbox{\C}^n$. 
It is totally complex if and only if the form 
$\frac{1}{2}\omega = \phi^*(\sum dx^i\wedge dy^i)$ 
is nondegenerate on $U$. We compute
\[ \frac{1}{2}\omega = \sum dx^i \wedge d \Re F_i = \sum 
(\Re \frac{\partial F_i}{\partial z^j}) dx^i \wedge dx^j - \sum 
(\Im  \frac{\partial F_i}{\partial z^j})dx^i \wedge du^j\, .\]  
{}From this it is easy to see that $\omega$ is nondegenerate if and 
only if the matrix $\Im  \partial F_i/\partial z^j$ is invertible,
i.e.\ if and only if $\phi$ is regular.   $\Box$ 

The following is a corollary of Lemma \ref{reg->Lemma}, 
Theorem \ref{ssconstrThm} and 
Theorem \ref{ssfundThm}. 
\begin{cor}\label{1Cor}
Any regular local holomorphic $1$-form $\phi$ on $\mbox{\C}^n$
defines a special symplectic manifold of complex dimension $n$. 
Conversely, any special symplectic manifold of complex dimension $n$
can be locally obtained in this way. 
\end{cor}
\begin{cor}\label{2Cor}
Any nondegenerate local holomorphic function on $\mbox{\C}^n$
defines a special K\"ahler manifold of complex dimension $n$. 
Conversely, any special K\"ahler manifold of complex dimension $n$
can be locally obtained in this way. 
\end{cor}
\noindent
{\bf Proof:}  A nondegenerate holomorphic function $F$ defines a 
regular and closed holomorphic $1$-form $dF$. The corresponding
holomorphic immersion $\phi = dF$ is totally complex 
and Lagrangian (by Lemma \ref{reg->Lemma}) and, by Proposition 
\ref{2->3Prop}, nondegenerate.  So it defines a 
special K\"ahler manifold by Theorem \ref{sKconstrThm}. The   
converse statement follows from Theorem~\ref{ssfundThm} 
and the fact that any holomorphic nondegenerate Lagrangian
immersion into $V$ is locally defined by a  
regular closed holomorphic
$1$-form (after choosing an appropriate isomorphism $V = T^*\mbox{\C}^n$). 
Notice that every  regular closed holomorphic
$1$-form on a simply connected domain is the differential of a nondegenerate 
holomorphic function. $\Box$ 

\section{Projective special geometry} 
\subsection{Conic and projective special manifolds}
We recall that a {\it local holomorphic} $\mbox{\C}^*$-{\it action} on a 
complex manifold $M$ is a holomorphic map 
\[ \mbox{\C}^*\times M \ni (\lambda , p) \mapsto 
\varphi_{\lambda}(p) \in M\]
defined on  an open neighbourhood $W$ of $\{ 1\} \times M$ such that
\begin{enumerate}
\item[(i)] $\varphi_1(p) = p$ for all $p\in M$ and
\item[(ii)] $\varphi_{\lambda}(\varphi_{\mu}(p)) = \varphi_{\lambda \mu}(p)$
if both sides are defined, i.e.\ if $(\lambda , \varphi_{\mu}(p)) \in W$
and $(\lambda \mu , p) \in W$. 
\end{enumerate} 
{}From this definition it follows that for every $p\in M$ there exist open
neighbourhoods $U_1$ of $1\in \mbox{\C}^*$ and $U_p$ of $p$ such that
$U_1 \times U_p \subset W$ and $\varphi_{\lambda}|U_p$ is a diffeomorphism
onto its image for all $\lambda \in U_1$. We will say that an equation
involving $\varphi_{\lambda}$ holds {\it locally} if it holds on any open set
$U \subset M$   on which $\varphi_{\lambda}$ is defined and on which it is 
a diffeomorphism onto its image. Of course, even if it is not explicitly 
mentioned, an equation involving $\varphi_{\lambda}$ is always meant  
to hold only locally. 

We use polar coordinates $(r,\theta )$ to parametrise 
$\mbox{\C}^* = \{ \lambda = re^{i\theta}| r,\theta \in\mbox{\R}, r>0\}$
and consider $\theta$ as a map from $\mbox{\C}^*$ to $\mbox{\R}/2\pi\mbox{\Z}$.

\begin{dof}\begin{enumerate}
\item[(i)] Let $(M,J,\nabla )$ be a complex manifold with a flat torsionfree
connection. It is called a {\it conic complex manifold}
if it admits a local holomorphic $\mbox{\C}^*$-action $\varphi_{\lambda}$ such
that locally 
$d\varphi_{\lambda}X = re^{\theta J}X = r(\cos \theta )X +
r(\sin \theta )JX$ for all $\n$-parallel vector fields $X$, 
where $\lambda = re^{i\theta}$.
\item[(ii)] A {\cmssll conic symplectic manifold} is a conic complex manifold
$(M,J,\n)$ together with a parallel symplectic form $\omega$. 
\item[(iii)] A conic symplectic manifold $(M,J,\n ,\omega )$ 
is called a {\cmssll conic K\"ahler manifold} if $\omega$ is $J$-invariant. 
\end{enumerate} 
\end{dof}
Notice that the condition $d\varphi_{\lambda}X = re^{\theta J}X$ for all
$\n$-parallel vector fields $X$ implies
that $\varphi_{\lambda}^*\n = \n^{\theta}$. 
\begin{prop} \label{conic->specialProp}
\begin{enumerate}
\item[(i)] Any conic complex manifold is a special complex manifold.
\item[(ii)] Any conic symplectic manifold is a special symplectic manifold. 
\item[(iii)] Any conic K\"ahler manifold is a special K\"ahler  manifold. 
\end{enumerate}
\end{prop}

\noindent
{\bf Proof:} Let $(M,J,\n )$ be a conic complex manifold and 
$\varphi_{\lambda}$ the corresponding local action. Since $d^{\nabla}J = 0$ 
is a local condition, it is sufficient to prove that any point $p\in M$
has an open neighbourhood $U$ such that $(U,J,\n )$ is a special
complex manifold. By Proposition \ref{prop1} it is sufficient to check that 
for any point $p\in M$ there
exist open neighbourhoods $U_1$ of $1\in \mbox{\R}/2\pi \mbox{\Z}$ 
and $U_p$ of $p$ such that $\n^{\theta}$ is a torsionfree
connection on $U_p$ for all $\theta \in U_1$.  From the definition
of local action it follows that for any $p\in M$ there exist 
open neighbourhoods $U_1$ of $1\in \mbox{\R}/2\pi \mbox{\Z}$ 
and $U_p$ of $p$ such that $\varphi_{\lambda}$  is defined on $U_p$ 
and $\varphi_{\lambda}|U_p$ is a diffeomorphism onto its image
for all $\lambda = e^{i\theta}$ with $\theta \in U_1$. Since 
$(M,J,\n )$ is a conic complex manifold we have  
$\n^{\theta} = \varphi_{\lambda}^*\n$ on $U_p$
for all $\theta \in U_1$.   Thus $\n^{\theta}$ is a 
torsionfree connection on $U_p$, proving (i). Statements (ii) and (iii) 
follow easily from (i). $\Box$

\begin{thm}\label{conicThm}
\begin{enumerate} 
\item[(i)] Let $(M,J,\n )$ be a complex manifold with a flat torsionfree
connection. Then $(M,J,\n )$ is a conic complex manifold if and only if
there exists a local holomorphic 
$\mbox{\C}^*$-action $\varphi_{\lambda}$ and  
for every $p\in M$  holomorphic functions $z^1,\ldots , z^n$
and $w_1,\ldots ,w_n$ defined near $p$  such that 
\begin{enumerate}
\item[(a)] 
$z^i\circ \varphi_{\lambda} = \lambda z^i$ and $w_j\circ \varphi_{\lambda} 
= \lambda w_j$ near $p$ and 
\item[(b)] $x^1 := \Re z^1, \ldots , x^n := \Re  z^n, y_1 := \Re w_1, 
\ldots , y_n := \Re w_n$
are affine local coordinates near $p$. 
\end{enumerate}
\item[(ii)] Let $(M,J,\n , \omega)$ be a complex manifold with 
a flat torsionfree connection and a parallel symplectic form. Then 
$(M,J,\n , \omega)$ is a conic symplectic manifold if and only if 
there exists a local holomorphic 
$\mbox{\C}^*$-action $\varphi_{\lambda}$ and  
for every $p\in M$  holomorphic functions $z^1,\ldots , z^n$
and $w_1,\ldots ,w_n$ defined near $p$  such that 
\begin{enumerate}
\item[(a)] 
$z^i\circ \varphi_{\lambda} = \lambda z^i$ and $w_j\circ \varphi_{\lambda} 
= \lambda w_j$ near $p$ and 
\item[(b)] $x^1 := \Re z^1, \ldots , x^n := \Re  z^n, y_1 := \Re w_1, 
\ldots , y_n := \Re w_n$
are affine local coordinates near $p$. 
\end{enumerate}
Moreover, if $(M,J,\n , \omega)$ is a conic (special) symplectic 
manifold then the local
holomorphic functions 
$z^i$ and $w_j$ can be chosen such that their real parts 
$x^i$ and $y_j$ form a real special coordinate system. 
\item[(iii)] Let $(M,J,\n , \omega)$ be a complex manifold with 
a flat torsionfree connection and a parallel $J$-invariant symplectic form.
Then 
$(M,J,\n , \omega)$ is a conic K\"ahler manifold if and only if 
there exists a local holomorphic  
$\mbox{\C}^*$-action $\varphi_{\lambda}$ and  
for every $p\in M$  holomorphic functions $z^1,\ldots , z^n$
and $w_1,\ldots ,w_n$ defined near $p$  such that 
\begin{enumerate}
\item[(a)] 
$z^i\circ \varphi_{\lambda} = \lambda z^i$ and $w_j\circ \varphi_{\lambda} 
= \lambda w_j$ near $p$ and 
\item[(b)] $x^1 := \Re z^1, \ldots , x^n := \Re  z^n, 
y_1 := \Re w_1, \ldots , y_n := \Re w_n$
are affine local coordinates near $p$. 
\end{enumerate}
Moreover, if $(M,J,\n , \omega)$ is a conic (special) K\"ahler
manifold then the local
holomorphic functions 
$z^i$ and $w_j$ can be chosen such that they form a conjugate pair of
special coordinates. 
\end{enumerate} 

\end{thm}

\noindent
{\bf Proof:} We prove only (i). Parts (ii) and (iii) are proven similarly. 
Let  $(M,J,\n )$ be a conic complex 
manifold and $x^1, \ldots , x^n, 
y_1, \ldots , y_n$ affine local coordinates on it. By Proposition 
\ref{conic->specialProp} and Theorem \ref{scoordThm} it is a special 
complex manifold and the affine local coordinates admit a 
holomorphic extension $z^1, \ldots , z^n, w_1, \ldots , w_n$. 
{}From $d\varphi_{\lambda}X = re^{\theta J}X$ for all $\n$-parallel
vector fields it follows that 
$z^i \circ \varphi_{\lambda} = \lambda z^i + c(\lambda)$, 
where $c : \mbox{\C}^* \rightarrow \mbox{\C}^n$ is a smooth
map. Since $\varphi_{\lambda}$ is a local action, the map $c$ must
satisfy the functional equation
\[ c(\lambda \mu) = \lambda c(\mu ) + c(\lambda )\]
for all $\lambda , \mu \in \mbox{\C}^*$ near $1 \in \mbox{\C}^*$ 
and $c(1) = 0$. It is easy to see that
this implies $c(\lambda ) = (1-\lambda )z_0$ for some constant
vector $z_0\in \mbox{\C}^n$. Up to adding (real) constants to the $x^i$,  
we can assume that the vector $z_0$ has purely imaginary components. 
Then changing the holomorphic extensions $z^i$ by adding purely imaginary
constants, we can arrange that
$c = z_0 = 0$ and hence that 
$z^i \circ \varphi_{\lambda} = \lambda z^i$. Similarly, we can show
that by adding constants one can arrange that $w_j \circ \varphi_{\lambda} 
= \lambda w_j$. This shows that a conic complex manifold
admits a local holomorphic $\mbox{\C}^*$-action and local 
holomorphic functions with the properties
(a) and (b).  Next we prove the converse statement of (i). 
So let  $\varphi_{\lambda}$ be a local holomorphic 
$\mbox{\C}^*$-action on $(M,J,\n )$ and  $z^1,\ldots , z^n, w_1,\ldots ,w_n$ 
local holomorphic functions satisfying (a) and (b).
{}From (a) and (b) it follows that $d\varphi_{\lambda}X = re^{\theta J}X$
for all $\n$-parallel vector fields $X$,  
by differentiation. This shows that $(M,J,\n )$ is a conic complex manifold. 
$\Box$ 

Next we are going to define the notion of projective special (complex,
symplectic or K\"ahler) manifold. These manifolds arise as orbit spaces
of conic special (complex, symplectic or K\"ahler) manifolds. 
Let $\varphi_{\lambda}$ be a local holomorphic $\mbox{\C}^*$-action
on a complex manifold $M$. To any point $p \in M$ we associate the
holomorphic curve $\varphi (p) : \lambda \mapsto \varphi_{\lambda}(p)$ in $M$
defined on an open neighbourhood of $1 \in \mbox{\C}^*$. If 
$\varphi_{\lambda}$ is the local $\mbox{\C}^*$-action associated to a
conic complex manifold then $\varphi (p)$ is an immersion and
${\cal D}_p := \varphi(p)_*T_1\mbox{\C}^* \subset T_pM$ defines an
integrable complex 1-dimensional holomorphic distribution on $M$. 
Its leaves are by definition the {\it orbits} of the local
$\mbox{\C}^*$-action $\varphi_{\lambda}$. We denote by
$\overline{M} = M/\mbox{\C}^*$ the set of orbits with the 
the quotient topology.  
$\overline{M}$  will be called the orbit space of $M$. 
If $M$ is a conic (complex, symplectic or K\"ahler) manifold and
the projection $M \ra \overline{M}$ is a holomorphic submersion onto a 
Hausdorff complex manifold, then $\overline{M}$ 
is called a {\it projective special} (complex, symplectic or K\"ahler) 
{\it manifold}. 

\subsection{Conic special coordinates} 
\begin{dof}\label{coniccoordDef} 
An affine local coordinate system 
$(x,y):=(x^1, \ldots , x^n, y_1, \ldots ,y_n)$
on a conic complex manifold $\,(M,J,\n )\,$ with corresponding 
local $\mbox{\C}^*$-action $\,\varphi_{\lambda}\,$ is called a 
{\cmssll conic affine local coordinate system} if it admits a holomorphic 
extension $\ (z,w):=$ $(z^1,\ldots,z^n,w_1,\ldots,w_n)$ such that 
locally $(z,w) \circ \varphi_{\lambda} = \lambda (z,w)$. Such a 
holomorphic extension is called a {\cmssll conic holomorphic extension}. 
\end{dof}

In view of Definition \ref{coniccoordDef} we will freely speak of
conic real special coordinate systems $(x,y)$ on conic symplectic 
manifolds and of conic conjugate pairs of special coordinates $(z,w)$
on conic K\"ahler manifolds. The following theorem is a corollary of 
Theorem \ref{conicThm}. 

\begin{thm} \label{ccoordThm} 
\begin{enumerate}
\item[(i)] Any conic complex manifold  
admits a conic local affine  
coordinate system near any point $p\in M$.  A conic local affine  
coordinate system is unique up to a linear transformation. 
\item[(ii)] Any conic symplectic manifold  
admits a conic real special coordinate system near any point $p\in M$.  
A conic real special coordinate system is unique up to a linear symplectic
transformation.  
\item[(iii)] Any conic K\"ahler manifold admits a conic conjugate
pair of special coordinates. A conic conjugate
pair of special coordinates is unique up to a (complex) linear symplectic
transformation.  
\end{enumerate} 
\end{thm}

\subsection{The extrinsic construction of conic and 
projective special manifolds} 
Let us consider the same fundamental data $V$, $\Omega$ and $\tau$
as in \ref{extrsmSec}. On $V$ we have the standard (global) holomorphic 
$\mbox{\C}^*$-action $\mbox{\C}^* \times V \ni (\lambda ,v) 
\mapsto \lambda v\in V$. 
A holomorphic immersion $\phi$ of a complex manifold $M$ into $V$
is called {\it conic} if for every point $p\in M$ and every neighbourhood
$U$ of $p$ there exist 
neighbourhoods $U_1$ of $1\in \mbox{\C}^*$ and $U_p$ of $p$ 
such that $\lambda \phi (U_p) \subset \phi (U)$ for all $\lambda \in U_1$.   
Notice that we do not require the image $\phi (M)$ to be a complex 
cone, i.e.\ (globally) invariant under the $\mbox{\C}^*$-action on $V$. 

\begin{thm}\label{csconstrThm}
Let $\phi$ be a conic totally complex holomorphic immersion of a complex
manifold $(M,J)$ into $V = T^* \mbox{\C}^n$, $n = \dim_{\mbox{\C}} M$,  
$\nabla$ the induced connection and $\omega = 2\phi^*(\sum dx^i\wedge dy_i)$ 
the induced symplectic form. Then the following hold:
\begin{enumerate}
\item[(i)] $(M,J,\nabla , \omega )$ is a conic  
symplectic manifold. 
\item[(ii)] The pull back via $\phi$ of the functions 
$(x^1 = \Re z^1, \ldots , x^n = \Re z^n, y_1 = \Re w_1, \ldots ,$\\ 
$y_n = \Re w_n)$ of $V$ defines a conic real special coordinate system 
around each point of $M$.  
\end{enumerate}
\end{thm}

\noindent
{\bf Proof:} Since $\phi$ is a conic holomorphic immersion, the 
holomorphic $\mbox{\C}^*$-action on $V$ induces a 
local holomorphic $\mbox{\C}^*$-action $\varphi_{\lambda}$ on $M$. 
One can easily check that $\varphi_{\lambda}$ defines on 
$(M,J,\nabla , \omega )$ the structure of a conic symplectic manifold
with conic real special coordinates $x^1\circ \phi , 
\ldots ,x^n\circ \phi , y_1\circ \phi , y_n\circ \phi$.   
$\Box$

\begin{thm}\label{cKconstrThm}
Let $\phi$ be a  conic holomorphic nondegenerate Lagrangian immersion of 
a complex manifold $(M,J)$ into $V$ inducing the K\"ahler metric $g$ on $M$.  
The immersion $\phi$ is totally complex and hence induces 
also the data $(\nabla , \omega )$ on $M$. Moreover, the following hold:
\begin{enumerate}
\item[(i)] $(M,J,\nabla , \omega )$ is a conic K\"ahler manifold.
\item[(ii)] $\omega$ coincides with the K\"ahler form of $g$, i.e.\ 
$\omega = g (\cdot , J\cdot )$.
\item[(iii)] The pull back via $\phi$ of the canonical  
coordinates  $(z^1, \ldots , z^n, w_1, \ldots , w_n)$ of $V$ defines 
a conic conjugate pair of special coordinates around each point of $M$. 
\end{enumerate}
\end{thm}

\noindent
{\bf Proof:} This follows from Theorem \ref{sKconstrThm}  and Theorem
\ref{csconstrThm}. $\Box$

Now we will show that any simply connected conic (complex, symplectic or
K\"ahler) manifold
arises by the construction of Theorem \ref{csconstrThm} or 
Theorem \ref{cKconstrThm}. 

\begin{thm}\label{csfundThm}
\begin{enumerate}
\item[(i)] Let  $(M,J,\nabla)$ be a simply connected conic 
complex manifold of complex dimension $n$. Then there exists a 
conic holomorphic totally complex  immersion $\phi : M \rightarrow 
V = T^* \mbox{\C}^n$ inducing the connection $\n$ on $M$. 
Moreover, $\phi$ is unique up to
a linear transformation of $V$ preserving the real structure  $\tau$. 
Here the real structure is considered as a (constant) field of antilinear 
involutions on the tangent spaces of $V$. Finally, $\omega = 2\phi^*
(\sum dx^i\wedge dy_i)$ is a $\n$-parallel symplectic structure 
defining on $(M,J,\n )$ the structure of conic symplectic manifold.  
\item[(ii)] Let  $(M,J,\nabla , \omega )$ be a simply connected conic  
symplectic manifold of complex dimension $n$. Then there exists a 
conic holomorphic totally complex  immersion $\phi : M \rightarrow 
V = T^* \mbox{\C}^n$ inducing the connection $\n$ and the
symplectic form $\omega$ on $M$. Moreover, $\phi$ is unique up to
a linear transformation of $V$ preserving the complex symplectic 
form $\Omega$ and the real structure  $\tau$.  
\item[(iii)]
Let  $(M,J,\nabla , \omega )$ be a simply connected conic 
K\"ahler manifold of complex dimension $n$ then there exists a 
conic holomorphic nondegenerate Lagrangian (and hence totally complex) 
immersion $\phi : M \rightarrow 
V = T^* \mbox{\C}^n$ inducing the K\"ahler metric $g$, the 
connection $\n$ and the symplectic form $\omega = 
2\phi^*(\sum dx^i\wedge dy_i) = g(\cdot , J\cdot )$ on $M$. Moreover 
$\phi$ is unique up to 
a linear transformation of $V$ preserving the complex symplectic 
form $\Omega$ and the real structure  $\tau$. Here the real structure
is considered as a field of antilinear involutions on the tangent spaces
of $V$.    
\end{enumerate}
\end{thm}

\noindent
{\bf Proof:} The proof is completely analogous to that of
Theorem \ref{ssfundThm}. To prove (ii), for instance, it is essentially 
sufficient to replace real special coordinates by conic real special 
coordinates in the proof of  Theorem \ref{ssfundThm} (ii). $\Box$

\label{cKfundThm}

We will call a holomorphic $1$-form $\sum F_idz^i$ on an 
open subset  $U \subset \mbox{\C}^n$ {\it conic} if 
the corresponding holomorphic immersion $U\ni z \mapsto \sum F_i(z)dz^i \in 
T^*_z\mbox{\C}^n \subset T^*\mbox{\C}^n = V$ is conic. This is the case
if and only if the functions $F_i$ are locally
homogeneous of degree one, i.e.\ if $F_i(\lambda z) 
= \lambda F_i(z)$ for all 
$z\in U$ and all $\lambda$ near $1\in \mbox{\C}^*$. 

A holomorphic function $F$ on $U$ is called {\it conic} if its differential
$dF$ is conic. This is the case if and only if $F$ is locally homogeneous
of degree 2, i.e.\ if $F(\lambda z) = \lambda^2 F(z)$ for all 
$z\in U$ and all $\lambda$ near $1\in \mbox{\C}^*$. 

We have the following analogues of Corollary \ref{1Cor} and 
Corollary \ref{2Cor}. 

\begin{cor}\label{c1Cor}
Any conic regular local holomorphic $1$-form $\phi$ on $\mbox{\C}^n$
defines a conic symplectic manifold of complex dimension $n$. 
Conversely, any conic symplectic manifold of complex dimension $n$
can be locally obtained in this way. 
\end{cor}
 
\begin{cor}\label{c2Cor}
Any conic nondegenerate local holomorphic function on $\mbox{\C}^n$
defines a conic K\"ahler manifold of complex dimension $n$. 
Conversely, any conic K\"ahler manifold of complex dimension $n$
can be locally obtained in this way. 
\end{cor}

\noindent
{\bf Remark 3:} Let  $\overline{M} = M/\mbox{\C}^*$ be  a projective 
special (complex, symplectic or K\"ahler) manifold, with
$M$ simply connected. Then the holomorphic immersion 
$\phi : M \rightarrow V$ constructed in Theorem \ref{csfundThm} 
induces a holomorphic immersion $\overline{\phi} : \overline{M} \ra P(V)$ 
into the complex projective space of complex dimension $2n-1$. 
The holomorphic immersion $\overline{\phi}$ is unique up to a projective 
transformation induced by a linear symplectic transformation of $V$ 
preserving the real structure $\tau$.  
To construct $\overline \phi$ it is sufficient to assume that 
$\overline M$ is simply connected.

\section{Geometric structures on the cotangent bundle of \ssm s}

In this section we prove that the cotangent bundle of a \ssm\
carries two canonical complex structures $J_1$, $J_2$. Moreover,  
if the (1,1)-part of the symplectic form $\o$ is
nondegenerate it also carries an almost hyper-Hermitian 
structure. This almost hyper-Hermitian 
structure is hyper-K\"ahler if and only if $\o^{11}$ is parallel. 
If the $(2,0)$-part of $\o$ is nondegenerate we obtain an
almost para-hypercomplex structure. It is para-hypercomplex 
if and only if $\o^{20}$ is parallel.  
This   generalises the known 
construction of a hyper-K\"ahler metric on the cotangent bundle of a
special K\"ahler manifold \cite{CFG,C2,F,H}.

Let $M$ be a manifold and denote by $N = T^*M$ its cotangent bundle. 
A connection $\n$ on $M$ defines a decomposition 
\be 
T_\x N = \ch^\n_\x  \oplus T^v_\x N \cong T_pM \oplus T^*_p M\ ,\ 
\x \in N\ ,\ p=\pi(\x)\ ,
\la{decomp}\ee
where $\pi : N=T^*M \ra M\,$, $T^v_\x N$ is the vertical subspace and  
$\ch^\n_\x $ is the horizontal subspace defined by the connection $\n$.
Here we have a natural identification of $T^v_\x N$ with $T^*_pM$
and an identification of $\ch^\n_\x$ with $T_p M$ defined by the 
projection $\pi$. If $M$ is a complex manifold with complex structure $J$, 
then $N$ carries a natural complex structure $J_N$. 
We note that the vertical subspace $T^v_\x N$ is $J_N$-invariant, 
but the horizontal subspace  $\ch^\n_\x$ is in general not.
We denote by $J^\n$ the almost complex structure on $N$ defined with respect
to the decomposition \re{decomp} by 
\be
J^\n = \pmatrix{J & 0 \cr
                0 & J^* }
\la{J1}\ee
In general $J^\n$ is not integrable. 

\bp \la{connProp}
Let $\n$ be a connection on a complex manifold $(M,J)$. The horizontal
distribution ${\ch}^\n \subset TN$ is $J_N$-invariant if and only if
there exists a torsionfree complex (i.e.\ $DJ = 0$) connection $D$
on $M$ such that the tensor field $A := \n - D$ satisfies the 
condition
\be A^\x_{X}\circ J =   A^\x_{JX}\quad \forall\ X\in TM\, ,\la{AJ}\ee 
where $A^\x_XY = \x (A_XY)$. 
\ep
For the proof we need two lemmas. The first one is well known.
\bl
Let $D$ and $\n$ be connections on a manifold $M$ and $A= \n - D$. Then the
corresponding horizontal distributions $\ch^D$ and $\ch^\n$ are related
by: 
\[ \ch^\n_\x = A^\x\ch^D_\x = \{ \wh{v} = v + A^\x_v| v\in \ch^D_\x \cong T_pM\}
\, ,\]
where $\x \in N = T^*M$ and $p = \pi (\x )$. 
\el 

\bl Let $\n$ be a torsionfree complex connection on a complex
manifold $(M,J)$. Then the horizontal
distribution ${\ch}^\n \subset TN$ is $J_N$-invariant and 
hence $J^\n = J_N$.
\la{DLemma}
\el

\noindent
{\bf Proof:} Let $(x^1,\ldots,x^n,y^1,\ldots,y^n,u_1,\ldots,u_n,v_1,\ldots,v_n)$ 
be the local coordinate system on $N=T^*M$ associated to a holomorphic local 
coordinate system $(z^1,\ldots,z^n)$ on $M$, i.e. $z^i = x^i + \sqrt{-1}y^i$
and $\o = \sum dx^i \wedge du_i + \sum dy^j \wedge dv_j $ is the canonical
symplectic structure on $N$. Note that 
$$ T^v N = 
\rm{span}\left\{\pp{}{u_1},\ldots\pp{}{u_n},\pp{}{v_1},\ldots,\pp{}{v_n}
                                                                 \right\}\ .
$$ 
We denote by $D$ the local connection on $M$ with horizontal space 
$$ \ch^D :=
 \rm{span}\left\{\pp{}{x^1},\ldots\pp{}{x^n},\pp{}{y^1},\ldots,\pp{}{y^n}
                                                                 \right\}\ .
$$ 
This connection is flat and torsionfree, with affine local coordinates
$x^1,\ldots,x^n,y^1,\ldots,y^n$. It is also complex because the complex 
structure $J$ is constant in these coordinates:
$$ 
J\pp{}{x^i} = \pp{}{y^i}\quad ,\quad  J\pp{}{y^i} = - \pp{}{x^i}\ .
$$
In terms of the induced coordinate system on $N$, $J_N$ is given by 
\bea
&& J_N\pp{}{x^i} =\ \ \pp{}{y^i}\quad ,\quad  J_N\pp{}{u_j} = -\pp{}{v_j}
\nonumber\\[6pt]
&& J_N\pp{}{y^i} = - \pp{}{x^i}\quad ,\quad   J_N\pp{}{v_j} = \pp{}{u_j} 
\nonumber\eea
This clearly shows that $J^D = J_N$.
Now let $\n$ be any torsionfree complex connection on $(M,J)$. This means that
the $(1,2)$ tensor $A=\n -D$ is symmetric and $J$-linear, i.e.
$$
A_X Y = A_Y X \quad,\quad [A_X\,,\,J] = 0\quad \forall\ X,Y\in TM\ .
$$ 
The latter equation can also be written in the form 
$\ J^* A_X^\x = A^{J^*\x}_X\ $ for all $\x\in T^*M$.  
We claim that this implies the $J_N$-invariance of $\, \ch^\n = A \ch^D$.
In fact we have 
$$\arr
J_N \wh v &=& J_N (v+ A^\x_v)\ =\ Jv + J^* A_v^\x\ =\ Jv + A_v^{J^*\x}\\[6pt]
 &=& Jv +  A^{J^*\x} v\ =\ Jv + J^* (A^\x)v \ =\  Jv + A^\x Jv \\[6pt] 
 &=& Jv + A^\x_{Jv} \ =\ \wh{Jv}\quad 
\forall v\in \ch^D_\x \cong T_pM\ ,\  p=\pi(\x)\in M\,. \quad\Box
\ea$$   

\noindent
{\bf Proof} (of Proposition  \ref{connProp}):
Let $D$ be a torsionfree complex connection on $M$ and $\n$ a connection
on $M$ such that $A=\n - D$ satisfies \re{AJ}. To prove that $\ch^\n$ is
$J_N$-invariant it suffices to check that 
$J_N\wh{v} = \wh{Jv}$ for all $v \in  \ch^D_\x \cong T_pM$. 
Using the identification \re{decomp} and the 
identity \re{AJ} we compute:
\[ J_N\wh{v} = J_N(v + A^\x_v) = Jv + J^*A^\x_v = 
Jv + A^\x_v \circ J = Jv + A^\x_{Jv} = \wh{Jv}\, .  \]  
Conversely, let $\n$ be a connection on $M$ such that $\ch^\n$ is 
$J_N$-invariant. From the integrability of $J$ it follows that there exists
a torsionfree complex connection $D$ on $M$. Now we check that
$J_N \ch^\n = \ch^\n $ implies \re{AJ}. For $\wh v = v+ A^\x_v  \in \ch^\n_\x$,
we have by Lemma~\ref{DLemma}:  $\ J_N \wh v = J v + J^* A^\x_v $. This shows
that $J_N \wh v \in  \ch^\n_\x$ if and only if  $J_N \wh v = \wh{Jv}$.
The latter equation is equivalent to  $ J^* A^\x_v = A^\x_{Jv} $, which
is precisely \re{AJ}. $\Box$ 

Now let $\omega$ be a field of nondegenerate bilinear forms on a manifold
$M$, considered
as a map $TM \ra T^*M$, and $\n$ a connection on $M$. Using the
identification \re{decomp} we define an almost complex structure
$J^\o$ on $N = T^*M$ by
\be \la{J2}  J^\o = \pmatrix{ 0 & -\o^{-1}\cr
                     \o & 0}
\ee

\bl \la{omegaLemma}
If $\n$ is flat and torsionfree and $\o$ is $\n$-parallel then
$J^\o$ is integrable.
\el

\noindent
{\bf Proof:} If we express $\,J^\o\,$ in terms of the canonical coordinates on 
$\,N = T^*M\,$ induced by local affine coordinates on $M$, then it  
has constant coefficients. This shows that $J^\o$ is integrable. $\Box$

\bt
Let $(M,J,\n ,\o )$ be a \ssm. Then the cotangent bundle $N = T^*M$ carries
two natural complex structures 
\[ J_1 = J^\n = \pmatrix{J & 0 \cr
                0 & J^* }\quad \mbox{and}\quad 
J_2 = J^\o = \pmatrix{ 0 & -\o^{-1}\cr
                     \o & 0}\quad .
\]  
The commutator and anticommutator of $J_1$ and $J_2$ are given by
\[ [J_1 , J_2 ] = 2J_1  \pmatrix{ 0 & -(\o^{-1})^{11}\cr
                     \o^{11} & 0}= -2\pmatrix{ 0 & -(\o^{-1})^{11}\cr
                     \o^{11} & 0} J_1\, , \] 
\[ \{ J_1, J_2\} = 2J_1  \pmatrix{ 0 & -(\o^{-1})'\cr
                     \o' & 0}= 2 \pmatrix{ 0 & -(\o^{-1})'\cr
                     \o' & 0}J_1\, ,\]
where $\o' = \omega^{20} + \o^{02}$. 
\et

\noindent
{\bf Proof:} The integrability of $J_2$ follows from Lemma 
\ref{omegaLemma}. To prove the integrability of $J_1$, by 
Proposition~\ref{connProp}, it is sufficient to check the identity
\re{AJ} for $A = \n{-}D = \frac{1}{2}J\n J$, 
where $D =\frac{1}{2}(\n + \n^{(J)})\,$ is the 
torsionfree complex connection of Proposition \ref{DProp}.
Using the fact that $\n J$ is symmetric we compute:
\[
 2A_X\circ J = J(\n_X J)\circ J = \n_X J = (\n_{\LARGE .} J)X = 
J(\n_{\LARGE .} J)JX = J(\n_{JX}J) = 2 A_{JX}\quad . \quad\Box 
\] 
 
\bt 
Let $(M,J,\n ,\o )$ be a \ssm.
\begin{enumerate}
\item[(i)] Assume that $\o^{11}$ is nondegenerate. Then the cotangent bundle 
$N = T^*M$ carries a canonical almost hyper-Hermitian structure 
$(J_1,J_2, J_3 = J_1J_2 = -J_2J_1, g_N)$ given by
\[ J_1 = J^\n \, , \quad J_2 = J^{\o^{11}}\, , \quad
 g_N = diag (g,g^{-1})\, ,\]
where $\ g = \o^{11} (J\cdot ,\cdot )$ is the K\"ahler metric
on $M$, see Proposition \ref{metricProp}. 
$J_1$ is the standard (integrable) complex
structure on the cotangent bundle of the complex manifold $(M,J)$.
The almost  hyper-Hermitian manifold
$\ (M,J_1,J_2, J_3,g_N)\ $ is hyper-Hermitian (i.e.\ the almost complex
structures $J_1,J_2, J_3$ are integrable) if and only if 
$\ \n \o^{11} =0$. In this case $\ (M,J_1,J_2, J_3,g_N)\,$ is a hyper-K\"ahler
manifold. 
\item[(ii)] Assume that $\o' = \o^{20} + \o^{02}$ is nondegenerate. 
Then the cotangent bundle 
$N = T^*M$ carries a canonical almost para-hypercomplex structure 
$(J_1,J_2)$, i.e.\ a commuting pair of almost complex structures, given by
\[ J_1 = J^\n \, , \quad J_2 = J^{\o'}\,   .\]
$J_1$ is again the standard (integrable) complex
structure and  $(J_1,J_2)$ is an (integrable) para-hypercomplex
structure (i.e.\ $J_1$ and $J_2$ is integrable) if and only if 
$\n \o^{20} = 0$. 
\end{enumerate} 
\et 

\noindent
Note that in the second case $J_3 =J_1J_2$ is not an almost complex structure
but an almost product structure, i.e.\ an involution. 

\noindent
{\bf Proof:}  Using the  identities
\[ J^* \circ \o^{11} = - \o^{11} \circ J\, ,\quad 
J^* \circ \o' = \o'\circ J\, ,\]
where the two-forms $\o^{11}$ and $\o'$ are considered as linear maps
$TM \rightarrow T^*M$, one can check that $J_1$ and $J_2$ are anticommuting
or commuting almost complex structures in case (i) and (ii) respectively.  
To check that $g_N$ is Hermitian with respect to the almost complex
structures $(J_1,J_2, J_3)$ in case (i) we compute $\o_{\alpha} :=
g_N \circ J_{\alpha}$ as follows:
\[ \o_1 = -\left( \sum \o_{ij}dq^i\wedge dq^j + \sum \o^{ij}dp_i\wedge dp_j 
\right) \, ,\]
where $\omega^{11} = \sum \o_{ij}(q)dq^i\wedge dq^j$ is the expression of the
symplectic form $\omega^{11}$ in $\n$-affine coordinates $q^i$ on $M$,
 $(\o^{ij}) = (\o_{ij})^{-1}$ and the $p_i$ are the  conjugate momenta
corresponding to the $q^i$.  
\[ \o_2 = \sum (J^*dq^j)\wedge dp_j \, , \quad 
\o_3 = \sum dq^j\wedge dp_j\, .\]
{}From these formulas we see that the $\o_{\alpha}$ are skew-symmetric
and therefore that the $J_{\alpha}$ are $g_N$-orthogonal. 
This shows that $(J_1,J_2, J_3,g_N)$ is an almost hyper-Hermitian structure.
The form $\o_3$ is closed. The form $\o_2$ is closed since $dJ^*\eta = 0$
for any parallel 1-form $\eta$.    
The form $\o_1$ is closed if and only if the coefficients $\o_{ij}$ are
constant, i.e. if and only if $\o^{11}$ is $\n$-parallel.  
If this is the case, 
then the almost hyper-Hermitian structure $(J_1,J_2,J_3,g_N)$  is 
hyper-K\"ahler, e.g.\ by Hitchin's Lemma.\par
Assume now that $J_2$ is integrable, i.e.\ 
the Nijenhuis tensor $N_{J_2}=0$. A direct calculation shows that
$$ J_2N_{J_2}(\partial_{q^i},\, \partial_{q^j}) = \sum_k(\rho_{jk,i} -
\rho_{ik,j})\partial_{p_k} \, ,
$$
where $\rho_{ij}(q)$ are the coefficients of $\rho = \o^{11}$ or $\o'$
in cases (i) or (ii) respectively. Notice that 
$\rho_{ik,j} - \rho_{jk,i}$ are the 
coefficients of the 2-form $d(i_{\partial_{q^k}}\rho) = 
L_{\partial_{q^k}} \rho$.  This shows that 
$N_{J_2}(\partial_{q^i},\, \partial_{q^j})=0$ implies that
Lie derivative of $\rho$ in the direction of any parallel vector field
on $M$ vanishes and hence that $\rho$ is parallel.   
$\Box$

\end{document}